# Calculating the sequences behind a hexagonal lattice based equal circle packing in the Euclidian plane

Jure Voglar [a] * and Aljoša Peperko [b, c]

[a] Department of Catalysis and Chemical Reaction Engineering, National Institute of Chemistry, Hajdrihova 19, 1001 Ljubljana, Slovenia, email: jure.voglar@ki.si

[b] Faculty of Mechanical Engineering, University of Ljubljana, Aškerčeva 6, 1000 Ljubljana, Slovenia, email: aljosa.peperko@fs.uni-lj.si

[c] Institute of Mathematics, Physics and Mechanics, Ljubljana, Jadranska 19, 1000 Ljubljana, Slovenia

* corresponding author

## Abstract

The article presents the mathematical sequences describing circle packing densities in four different geometric configurations involving a hexagonal lattice based equal circle packing in the Euclidian plane. The calculated sequences take form of either polynomials or rational functions. If the circle packing area is limited with a circle, the packing densities tend to decrease with increasing number of the packed circles and converge to values lower than $\pi/(2\sqrt{3})$. In cases with packing areas limited by equilateral triangles or equilateral hexagons the packing densities tend to increase with increasing number of the packed circles and converge to $\pi/(2\sqrt{3})$. The equilateral hexagons are shown to be the preferred equal circle packing surface areas with practical applications searching for high equal circle packing densities, since the packing densities with circle packing inside equilateral hexagons converge faster to $\pi/(2\sqrt{3})$ than in the case of equilateral triangle packing surface areas.

## Keywords

circle packing; Euclidian plane; hexagonal lattice; mathematical sequences; practical applications

# 1. Introduction

Circular packing has a wide range of potential applications [1]. The examples of applications span from design of conductor cables [2], logistics sector [3–8] and robotics [9]. Specifically, the equal circles packing in confined surface area geometries has been known to have the highest use value in the industry, since multiple mass-produced products or refined materials have the same shape (e.g. cylinder or disc) and size (e.g. paper rolls, hot rolled steel coils, rolled textile, wire and cable rolls) and have to be efficiently stored and transported usually in ships by large distances. Also, appropriate version of equal circle packing can serve as a guide in design of monolith catalysts used in chemical industry [10] and/or environmental catalysis [11].

In 1773, J.L. Lagrange proved that the highest-density lattice packing of circles in the two-dimensional Euclidean plane is the hexagonal packing arrangement [12], where the centers of circles are arranged in a hexagonal lattice and each circle is surrounded by six other circles. It is easy to see that this packing density equals to $\frac{\pi}{2\sqrt{3}} \approx 0.9069$. In 1890, A. Thue attempted to prove that this density is optimal among all packings, but his proof was incomplete. In 1942, L.F. Tóth published the first rigorous proof [12,13].

Recently, various different numerical algorithms were obtained and tested to achieve high circle packing densities in various different geometrically restricted surface areas [14–19]. Our investigation is focused on a hexagonal lattice based equal circle packing by calculating the packing density sequences in chosen geometric settings. The determinations of packing density can be particularly useful to serve as a decision-making tool, when a designer needs to arrange equally sized objects into a confined space with a sufficiently high packing density.

# 2. Methods

Our work started with derivations for the first 6 cases (with increasing of the number of packed circles) of a hexagonal lattice based equal circle packing arrangements. The derivations involved visual observations of different packing configurations and extractions of the number of packed circles, geometric properties, relations between the packed circles and the bounding geometry and calculations of the circle packing density. Four different cases were considered: (a) circle packing forming tringle inside a bounding circle with center in the center of the triangle (Figure 1 a), (b) circle packing forming tringle inside a bounding equilateral triangle (Figure 1 b), (c) circle packing forming hexagon inside a bounding circle with center in the center of the hexagon (Figure 1 c) and (d) circle packing forming hexagon inside a bounding equilateral hexagon (Figure 1 d).

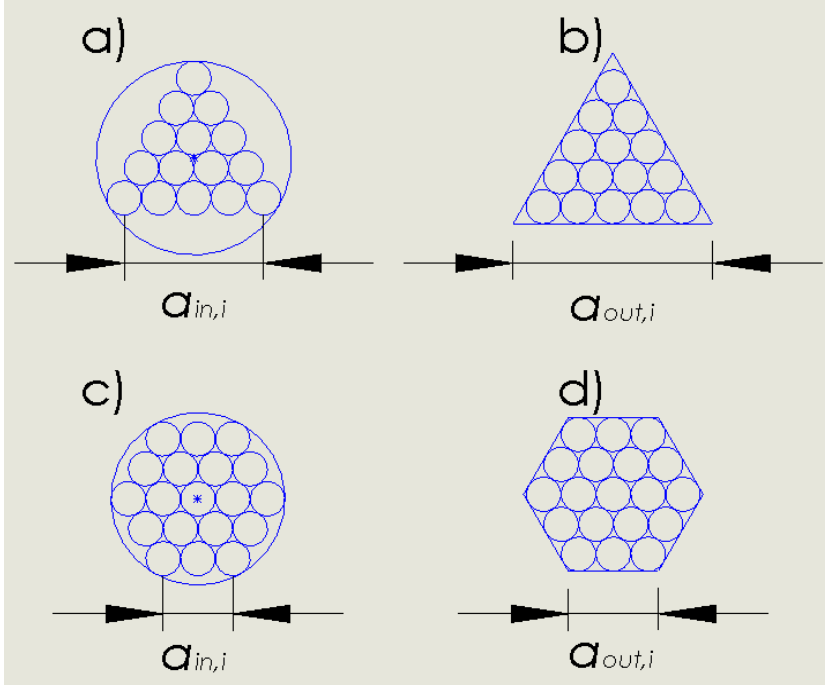

Figure 1. Different hexagonal lattice based equal circle packing configurations: (a) and (b) at *i* = 5, (c) and (d) at *i* = 2.

After the initial derivations and calculations, we started to investigate the trends in the observed quantities, especially the number of packed circles and the geometric relations between the packed circles and the bounding geometry.

## 3. Results and discussion

The calculated sequences depend on the index number *i*, indexing the number of the packing arrangement with the packed circles of radii *r*.

### 3.1. Sequences of the (a) circle packing

The number of the packed circles $N_i$ follows the triangular numbers sequence [20]:

$$N_i = \frac{i(i+1)}{2} = \frac{i^2+i}{2}.$$

The length of the side of the equilateral triangle of the packed circles $a_{in,i}$ follows the even number sequences:

$$a_{in,i} = 2(i-1)r \ .$$

The radius of the surrounding circle $R_i$ is equal to:

$$R_i = \frac{a_{in,i}}{\sqrt{3}} + r = \frac{2(i-1)r}{\sqrt{3}} + r = r\left[\frac{2(i-1)}{\sqrt{3}} + 1\right] = r\left[\frac{2i-2+\sqrt{3}}{\sqrt{3}}\right].$$

The ratio of the radii $R_i/r$ is therefore:

$$\frac{R_i}{r} = \frac{2(i-1)}{\sqrt{3}} + 1 = \frac{2i-2+\sqrt{3}}{\sqrt{3}}.$$

The packing density $\rho_{p,a,i}$ equals:

$$\rho_{p,a,i} = \frac{N_i \pi r^2}{\pi R_i^2} = N_i \left(\frac{r}{R_i}\right)^2 = \frac{N_i}{\left(\frac{R_i}{r}\right)^2} = \frac{\frac{i^2+i}{2}}{\left(\frac{2i-2+\sqrt{3}}{\sqrt{3}}\right)^2} = \frac{\frac{i^2+i}{2}}{\left(\frac{2i-(2-\sqrt{3})}{\sqrt{3}}\right)^2} = \frac{\frac{i^2+i}{2}}{\frac{4i^2-4i(2-\sqrt{3})+(2-\sqrt{3})^2}{3}} =$$

$$\frac{3i^2+3i}{8i^2+(8\sqrt{3}-16)i+(14-8\sqrt{3})} = \frac{3}{8} + \frac{(9-3\sqrt{3})i-\left(\frac{21}{4}-3\sqrt{3}\right)}{8i^2+(8\sqrt{3}-16)i+(14-8\sqrt{3})}.$$

The packing density's limit as the number of circles tend to infinity is equal to:

$$\lim_{i\to\infty} \rho_{p,a,i} = \frac{3}{8} = 0.375.$$

The sequences describing the number of packed circles, length of sides of the equilateral triangles of the packed circles, radii of the outer circle and the ratios of the radii take a form of different polynomials, while the packing density takes a form of a rational function.

### 3.2. Sequences of the (b) circle packing

The number of the packed circles $N_i$ follows the triangular numbers sequence as in the case (a).

The length of the side of the outer equilateral triangle $a_{out,i}$ is:

$$a_{out,i} = 2(i-1)r + 2\sqrt{3}r = \left(2i-1+2\sqrt{3}\right)r,$$

where $i > 1$.

The surface area of the outer equilateral triangle $A_i$ equals:

$$A_i = \frac{\sqrt{3}}{4} a_{out,i}^2 = \frac{\sqrt{3}}{4}\left[\left(2i-1+2\sqrt{3}\right)r\right]^2 = \frac{\sqrt{3}}{4}\left[\left(2i-(1-2\sqrt{3})\right)r\right]^2 = \frac{\sqrt{3}}{4}\left[4i^2 - \left(4-8\sqrt{3}\right)i+(13-4\sqrt{3})\right]r^2 = \left[\sqrt{3}i^2-\left(\sqrt{3}-6\right)i+\left(\frac{13\sqrt{3}}{4}-3\right)\right]r^2.$$

The packing density $\rho_{p,b,i}$ inside the outer equilateral triangle is:

$$\rho_{p,b,i} = \frac{N_i \pi r^2}{A_i} = \frac{\left(\frac{i^2+i}{2}\right)\pi r^2}{\left[\sqrt{3}i^2-(\sqrt{3}-6)i+(\frac{13\sqrt{3}}{4}-3)\right]r^2} = \frac{\pi i^2+\pi i}{2\left[\sqrt{3}i^2-(\sqrt{3}-6)i+(\frac{13\sqrt{3}}{4}-3)\right]} = \frac{\pi}{2\sqrt{3}} + \frac{\pi\left[\left(\frac{2}{\sqrt{3}}-2\right)i-\frac{13}{4\sqrt{3}}+1\right]}{2i^2-(2-4\sqrt{3})i+(\frac{13}{2}-2\sqrt{3})}.$$

The packing density's limit at infinite number of circles inside the outer equilateral triangle is:

$$\lim_{i\to\infty} \rho_{p,b,i} = \frac{\pi}{2\sqrt{3}} \approx 0.9069.$$

Therefore, the packing density inside the outer equilateral triangle converges as in Thule's theorem [12] to the value of the highest-density lattice packing of equal circles in the Euclidian plane.

The sequences describing the lengths of the sides of the outer equilateral triangles and the surface areas of the outer equilateral triangles take a form of different polynomials, while the packing density takes a form of a rational function.

The Table 1 illustrates that the packing density falls in case (a) and raises in case (b) with increase of the number of packed circles.

Table 1. Summary of the first values of the (a) and (b) packing arrangements.

| $i$ | $N_i$ | $R_i/r$ | $\rho_{p,a,i}$ | $\rho_{p,b,i}$ |
|---|---|---|---|---|
| 1 | 1 | 1 | 1 | / |
| 2 | 3 | 2.154701 | 0.646171 | 0.520899741 |
| 3 | 6 | 3.309401 | 0.547838 | 0.607629359 |
| 4 | 10 | 4.464102 | 0.501801 | 0.662590706 |
| 5 | 15 | 5.618802 | 0.475121 | 0.700516766 |
| 6 | 21 | 6.773503 | 0.457712 | 0.728258568 |
| 7 | 28 | 7.928203 | 0.44546 | 0.749430204 |
| 8 | 36 | 9.082904 | 0.436368 | 0.766117485 |
| 9 | 45 | 10.2376 | 0.429354 | 0.779608173 |
| 10 | 55 | 11.3923 | 0.423779 | 0.790740196 |
| 11 | 66 | 12.54701 | 0.419241 | 0.800082236 |
| 12 | 78 | 13.70171 | 0.415475 | 0.808033801 |
| 13 | 91 | 14.85641 | 0.4123 | 0.814883767 |
| 14 | 105 | 16.01111 | 0.409587 | 0.820846187 |
| 15 | 120 | 17.16581 | 0.407242 | 0.826083037 |
| 16 | 136 | 18.32051 | 0.405195 | 0.830719148 |
| 17 | 153 | 19.47521 | 0.403392 | 0.834852268 |
| 18 | 171 | 20.62991 | 0.401792 | 0.838560035 |
| 19 | 190 | 21.78461 | 0.400363 | 0.841904891 |
| 20 | 210 | 22.93931 | 0.399079 | 0.844937627 |
| 21 | 231 | 24.09401 | 0.397918 | 0.84769998 |
| 22 | 253 | 25.24871 | 0.396864 | 0.850226562 |
| 23 | 276 | 26.40341 | 0.395903 | 0.852546323 |
| 24 | 300 | 27.55811 | 0.395023 | 0.854683656 |
| 25 | 325 | 28.71281 | 0.394214 | 0.856659266 |

The circle packing (b) could be for instance useful for transport of fragile cylinder-shaped objects inside e.g. boats.

### 3.3. Sequences of the (c) circle packing

The number $N_i$ of the packed circles equals:

$$N_i = N_{i-1} + 6i \,,$$

$$N_i = N_0 + \frac{i(i+1)}{2} \cdot 6 = N_0 + 3 \cdot i(i+1) = N_0 + 3(i^2 + i),$$

where $N_0 = 1$ and so:

$$N_i = 1 + 3(i^2 + i) = 1 + 3i^2 + 3i = 3i^2 + 3i + 1.$$

The length of the side of the equilateral hexagon of circles $a_{in,i}$ is another even number sequence:

$$a_{in,i} = 2ir.$$

The ratio of the radii is a sequence of odd numbers:

$$\frac{R_i}{r} = 2i + 1.$$

The packing density $\rho_{p,c,i}$ equals:

$$\rho_{p,c,i} = \frac{N_i \pi r^2}{\pi R_i^2} = N_i \left(\frac{r}{R_i}\right)^2 = \frac{N_i}{\left(\frac{R_i}{r}\right)^2} = \frac{3i^2 + 3i + 1}{(2i+1)^2} = \frac{3i^2 + 3i + 1}{4i^2 + 4i + 1} = \frac{3}{4} + \frac{1}{16i^2 + 16i + 4}.$$

The packing density's limit as the number of circles tends to infinity is:

$$\lim_{i \to \infty} \rho_{p,c,i} = \frac{3}{4} = 0.75 \,.$$

The sequences describing the number of packed circles, lengths of the sides of the equilateral hexagons of the packed circles and ratios of the radii take a form of different polynomials, while the packing density takes form of a rational function.

### 3.4. Sequences of the (d) circle packing

The number of the packed circles $N_i$ is the same as in the case (c).

The length of side of the outer equilateral hexagon $a_{out,i}$ for $i > 0$ is:

$$a_{out,i} = 2ir + \frac{2}{\sqrt{3}}r = \left(2i + \frac{2}{\sqrt{3}}\right)r \,.$$

The surface area of the outer equilateral hexagon $A_i$ is then:

$$A_i = \frac{3\sqrt{3}}{2} a_{out,i}^2 = \frac{3\sqrt{3}}{2}\left[\left(2i + \frac{2}{\sqrt{3}}\right)r\right]^2 = \frac{3\sqrt{3}}{2}(4i^2 + \frac{8}{\sqrt{3}}i + \frac{4}{3})r^2 = (6\sqrt{3}i^2 + 12i + \frac{4\sqrt{3}}{2})r^2.$$

The packing density $\rho_{p,d,i}$ inside the outer equilateral hexagon equals:

$$\rho_{p,d,i} = \frac{N_i \pi r^2}{A_i} = \frac{(3i^2+3i+1)\pi r^2}{(6\sqrt{3}i^2+12i+\frac{4\sqrt{3}}{2})r^2} = \frac{3\pi i^2+3\pi i+\pi}{6\sqrt{3}i^2+12i+\frac{4\sqrt{3}}{2}} = = \frac{\pi}{2\sqrt{3}} + \frac{\pi(1-\frac{2}{\sqrt{3}})i}{2\sqrt{3}i^2+4i+\frac{2}{\sqrt{3}}}.$$

The packing density's limit as the number of circles tends to infinity is:

$$\lim_{i \to \infty} \rho_{p,d,i} = \frac{\pi}{2\sqrt{3}} \approx 0.9069.$$

This limit is the same as in the case (b) and describes the convergence to the value of the highest-density lattice packing of equal circles in the Euclidian plane.

The sequences describing the lengths of the outer equilateral hexagons and the surface areas of the outer equilateral hexagons take a form of different polynomials, while the packing density takes a form of a rational function.

The Table 2 illustrates that the packing density falls in case (c) and rises in case (d) with increase of the number of packed circles. By comparing the number of packed circles with the Table 1, we can clearly observe much faster increase of $N_i$ with the circles arranged in hexagons. For instance, the 25th element of the Table 1 and Table 2 have the values of 325 and 1801, respectively.

Table 2. Summary of the first values of the (c) and (d) packing arrangements.

| $i$ | $N_i$ | $R_i/r$ | $\rho_{p.c.i}$ | $\rho_{p.d.i}$ |
|---|---|---|---|---|
| 0 | 1 | 1 | 1 | / |
| 1 | 7 | 3 | 0.777778 | 0.850511 |
| 2 | 19 | 5 | 0.76 | 0.864659 |
| 3 | 37 | 7 | 0.755102 | 0.874011 |
| 4 | 61 | 9 | 0.753086 | 0.880115 |
| 5 | 91 | 11 | 0.752066 | 0.884349 |
| 6 | 127 | 13 | 0.751479 | 0.887442 |
| 7 | 169 | 15 | 0.751111 | 0.889795 |
| 8 | 217 | 17 | 0.750865 | 0.891644 |
| 9 | 271 | 19 | 0.750693 | 0.893134 |
| 10 | 331 | 21 | 0.750567 | 0.89436 |
| 11 | 397 | 23 | 0.750473 | 0.895386 |
| 12 | 469 | 25 | 0.7504 | 0.896257 |
| 13 | 547 | 27 | 0.750343 | 0.897006 |
| 14 | 631 | 29 | 0.750297 | 0.897656 |
| 15 | 721 | 31 | 0.75026 | 0.898227 |
| 16 | 817 | 33 | 0.75023 | 0.898731 |
| 17 | 919 | 35 | 0.750204 | 0.89918 |
| 18 | 1027 | 37 | 0.750183 | 0.899582 |
| 19 | 1141 | 39 | 0.750164 | 0.899945 |
| 20 | 1261 | 41 | 0.750149 | 0.900273 |
| 21 | 1387 | 43 | 0.750135 | 0.900572 |
| 22 | 1519 | 45 | 0.750123 | 0.900844 |
| 23 | 1657 | 47 | 0.750113 | 0.901095 |
| 24 | 1801 | 49 | 0.750104 | 0.901325 |

## 4. Conclusions

In the article we calculated the mathematical sequences behind a hexagonal lattice based equal circle packing in the Euclidian plane. Four different circle packing arrangements were studied and analyzed. The results revealed the following conclusions:

- The circles forming equilateral triangles (cases a and b) exhibited lower packing densities compared to the cases when they form equilateral hexagons (cases c and d).
- When the packing areas were limited by circles (cases a and c) the packing density showed monotonically decreasing trend, while in cases the areas were limited by equilateral polygons (triangles and hexagons, cases b and d) the packing density monotonically increased (with increasing value of $i$).
- The packing density's limit values in cases of circle bounded areas (cases a and c) were lower (3/8 and 3/4) than the value of the highest-density lattice packing of $\pi/(2\sqrt{3})$. The reason for this behavior is the fact that some proportion of the bounding circle's surface area is always unoccupied by the packed circles.

- The packing density's limit values in cases of equilateral polygons (triangles and hexagons, cases b and d) bounded areas were both equal to the value of the highest-density lattice packing of $\pi/(2\sqrt{3})$. In this case the packing limit takes the highest possible value of hexagonal equal circle packing, since the unoccupied surface area near the edges of the bounding polygons gets more and more covered by the packed circles with increasing of the number $i$.

- From the presented results we can conclude that in case of a hexagonal lattice based equal circle packing the highest packing densities can be achieved by packing them in areas forming equilateral triangles and/or equilateral hexagons (cases b and d). The packing densities with circle packing inside equilateral hexagons (case d) converge faster to the $\pi/(2\sqrt{3})$ limit and should thus be the preferred option in multiple practical applications. This notion could be an important piece of knowledge for advancements and optimizations in design of many industrial and logistics applications e.g. warehouses, shipping containers etc.

## Acknowledgements

J.V. thanks the support from the Slovenian Research and Innovation Agency (ARIS) through core funding P2-0152.

A.P. acknowledges the partial support of the Slovenian Research and Innovation Agency (grant P1-0222).